\documentclass[12pt]{article}
\usepackage{amssymb}
\usepackage{amsmath}
\usepackage{latexsym}
\usepackage{amsthm}
\usepackage{mathptmx}

\makeatletter
\def\thmhead@plain#1#2#3{%
  \thmname{#1}\thmnumber{\@ifnotempty{#1}{ }#2}%
  \thmnote{ \the\thm@notefont(#3)}}
\let\thmhead\thmhead@plain
\def\swappedhead#1#2#3{%
  \thmnumber{#2}\thmname{\@ifnotempty{#2}{. }#1}%
  \thmnote{ \the\thm@notefont(#3)}}
\makeatother

\theoremstyle{definition} 

 \newtheorem{definition}{Definition}[section]

\theoremstyle{plain}      

 \newtheorem{proposition}[definition]{Proposition}
 \newtheorem{theorem}[definition]{Theorem}
 \newtheorem{corollary}[definition]{Corollary}
 \newtheorem{lemma}[definition]{Lemma}
 
 \def \dem{\noindent{\sc Proof.~}}
\def \findem {\hfill{\hbox {\vrule\vbox{\hrule width 6pt\vskip 6pt\hrule}\vrule}}}

\def\OC{{\mathcal{O}}}
\def\FC{{\mathcal{F}}}
\def\g{{\mathfrak{g}}}
\def\d{{\mathfrak{d}}}
\def\m{{\mathfrak{m}}}

\def\NM{{\mathbb{N}}}

\def \nn{\{1,\dots,n\}}
\def \pp{\{1,\dots,p\}}
\def \uk{\underline{k}}

\def \deg {\mathop{\hbox{\rm deg}}\nolimits}
\def \Card {\mathop{\hbox{\rm card}}\nolimits}
\def \supp {\mathop{\hbox{\rm supp}}\nolimits}

\def \Ad {\mathop{\hbox{\rm Ad}}\nolimits}
\def \Im {\mathop{\hbox{\rm Im}}\nolimits}

\def \mod {\mathop{\hbox{\rm mod}}\nolimits}
\def \ad {\mathop{\hbox{\rm ad}}\nolimits}
\def \op {{\scriptstyle{\rm op}}}
\def \tot {{\scriptstyle{\rm tot}}}
\def \Id {\mathop{\hbox{\rm id}}\nolimits}
\def \Oplus {\mathop{\oplus}\limits}
\def \lar {\mathop{\leftarrow}\limits}
\def \To {\mathop{\longrightarrow}\limits}
\def \Otimes {\mathop{\otimes}}
\def \Fnp {\mathop{\hbox{\rm Free}}\nolimits_{(n,p)}}
\def \F {\mathop{\hbox{\rm Free}}\nolimits}

\def \ve{\wedge}

\begin{document}
\title{A $\hbar$-adic valuation property \\ of universal $R$-matrices}

\author{Benjamin Enriquez${}^\dag$ and Gilles Halbout${}^\ddag$\\
{\small{Institut de Recherche Math\'ematique Avanc\'ee de Strasbourg}}\cr
{\small{UMR 7501 de l'Universit\'e Louis Pasteur et du CNRS}}\cr
{\small{7, rue R. Descartes F-67084 Strasbourg }}\cr
{\small{${}^\dag$ \! e-mail:\,\texttt{enriquez@math.u-strasbg.fr}}}\cr
{\small{${}^\dag$ \! e-mail:\,\texttt{halbout@math.u-strasbg.fr}}}}
\markboth{Benjamin Enriquez and Gilles Halbout}
{A $\hbar$-adic valuation property of universal $R$-matrices}

\maketitle

\abstract{We prove that if $U_\hbar(\g)$ is a quasitriangular QUE
algebra with universal $R$-matrix $R$,
and $\OC_\hbar(G^*)$ is the quantized function algebra sitting inside
$U_\hbar(\g)$, then $\hbar\log(R)$ belongs to the tensor square $\OC_\hbar(G^*)
\bar{\otimes}\OC_\hbar(G^*)$. This gives another proof of the results of
Gavarini and Halbout, saying that $R$ normalizes
$\OC_\hbar(G^*) \bar{\otimes} \OC_\hbar(G^*)$ 
and therefore induces a braiding of the formal group $G^*$
(in the sense of Weinstein and Xu, or
Reshetikhin).}

\vskip20pt

\centerline {\bf \S\; 0 \ Introduction }
\label{introduction}

\vskip20pt

Let $(\g,r)$ be a finite-dimensional quasitriangular Lie bialgebra over a field
$k$ of characteristic $0$, and let $(U_\hbar(\g),R)$ be a 
quasitriangular quantization 
of $(\g,r)$. Recall that this means that (see \cite{Dr}):

\noindent$(1)$ $(\g,[-,-],\delta)$ is a Lie bialgebra, $r\in \g \otimes \g$ is
a solution of the classical Yang-Baxter equation, and the cobracket 
$\delta$~: $\g \to \ve^2(\g)$ of $\g$ is given by
$\delta(x)=[x \otimes 1 + 1 \otimes x,r]$ for $x \in \g$;

\noindent$(2)$ $(U_\hbar(\g),m,\Delta)$ is a quantized universal enveloping 
(QUE) algebra
($m$ is the product of $U_\hbar(\g)$, $\Delta$ is its coproduct), 
$R \in U_\hbar(\g)^{\widehat{\otimes} 2}$ satisfies the quasitriangular identities:
$$(\Delta \otimes \Id)(R)=R^{13}R^{23},\hskip1.5cm (\Id \otimes
\Delta)(R)=R^{13}R^{12},$$
$$(\Id \otimes \varepsilon)(R)=(\varepsilon \otimes \Id)(R)=1,\hskip1cm
\Delta^{\op}=\Ad(R)\circ \Delta,$$
where $\Ad(R)(x)=RxR^{-1}$ for $x \in U_\hbar(\g)^{\widehat{\otimes} 2}$, and
\begin{equation}
\left({{1}\over{\hbar}}(R-1)~\mod \hbar\right)=r~; 
\label{R:r}
\end{equation}
here $\widehat{\otimes}$ denotes the $\hbar$-adically completed tensor product,
the map $x \mapsto (x \mod \hbar)$ is the canonical projection
$$U_\hbar(\g)^{\widehat{\otimes} 2}\to U_\hbar(\g)^{\widehat{\otimes} 2}
\Otimes\limits_{k[[\hbar]]} k=U(\g)^{{\otimes} 2},$$ 
and $r$ is viewed as an
element of $U(\g)^{{\otimes} 2}$.
For $n \geq 0$, let us denote by $\delta_n$~: $U_\hbar(\g) \to
U_\hbar(\g)^{\widehat{\otimes} n}$ the map 
$$\delta_n =(\Id - \eta \circ \varepsilon)^{\otimes n} \circ \Delta^{(n)},$$
where  $\Delta^{(n)}$~: $U_\hbar(\g) \to U_\hbar(\g)^{\widehat{\otimes} n}$
is the $n$-th fold coproduct and $\eta$, $\varepsilon$ are the unit
and counit maps of
$U_\hbar(\g)$. Set
$$\OC_\hbar(G^*)=\{x \in U_\hbar(\g)~|~\forall n \geq 0,~\delta_n(x) \in \hbar^n 
~U_\hbar(\g)^{\widehat{\otimes} n}\}.$$
Then a classical result (see [Dr, Ga]) says that $\OC_\hbar(G^*)$
is a quantization of the Hopf-Poisson algebra $\OC(G^*)$ of functions
on the formal group corresponding to 
$\g^*$
(so $\OC(G^*)=U(\g^*)^*$).

\smallskip

According to $(\ref{R:r})$, we have $R \in 1+\hbar U_\hbar(\g)^{\widehat{\otimes} 2}$.
So $\log(R)=\sum\limits_{n \geq 1} {(-1)}^n {{(R-1)^n}\over{n}}$ is
a well-defined element of $\hbar 
U_\hbar(\g)^{\widehat{\otimes} 2}$.

\begin{theorem}
\label{thm:main}
Set $\rho={{1}\over{\hbar}} \log(R)$. Then
$\rho$ belongs to
$\OC_\hbar(G^*)^{\bar{\otimes}2}$ (this is a subalgebra of
$U_\hbar(\g)^{\widehat{\otimes} 2}$). If $\m_\hbar$ is the augmentation ideal
of $\OC_\hbar(G^*)$, we even have $\rho \in (\m_\hbar)^{\bar{\otimes}2}$.
\end{theorem}

Here $\bar{\otimes}$ denotes the
completed tensor product of formal series algebras.
As a corollary, we obtain a result of \cite{GH}.

\begin{corollary} (\cite{GH})
\label{cor:GH}
The $R$-matrix $R$ normalizes $\OC_\hbar(G^*)^{\bar{\otimes}2}$, in other words,
$\Ad(R)$~:
$U_\hbar(\g)^{\widehat{\otimes} 2} \to U_\hbar(\g)^{\widehat{\otimes} 2}$
restricts to an automorphism of 
$\OC_\hbar(G^*)^{\bar{\otimes}2}$.
\end{corollary}

Gavarini and one of us (see \cite{GH}) derive from this result that 
$\Ad(R)_{|\hbar=0}$,
the reduction $\mod \hbar$ of $\Ad(R)$ is an 
automorphism of $\OC_\hbar(G^*)\bar{\otimes}\OC_\hbar(G^*)$,
satisfying the braiding identities of \cite{WX} or \cite{Re}.
In a forthcoming paper, we plan to prove that this braiding coincides with
the braiding defined by 
Weinstein and Xu in \cite{WX}.

\medskip

This paper is organized as follows:
we prove Theorem \ref{thm:main} in Section 2.b.
This proof uses a combinatorial result on universal Lie algebras $F_{(n,p)}$,
which is stated and proved in Section 1.
In Section 3, we prove Corollary \ref{cor:GH} using
Theorem \ref{thm:main}.

\pagebreak

\centerline {\bf \S\; 1 \ A theorem on the universal Lie algebra $F_{(n,p)}$}
\label{part1}

\stepcounter{section}

\vskip20pt

In this section, we introduce a Lie algebra $F_{(n,p)}$
(Section 1.a). This Lie algebra is universal for the following
situation: $A$ is an algebra, $\rho \in A^{\otimes 2}$, and we consider
elements $\rho^{i,j} \in A^{\otimes (n+p)}$, $i \in\nn$, $j\in\{n+1,\dots,n+p\}$.
We construct elements $\delta_{(n,p)}$ in a completion of $F_{(n,p)}$, using
the Campbell-Baker-Hausdorff series (Section 1.b). The main result is Theorem
\ref{valuation} on the valuation (for the total degree) of $\delta_{(n,p)}$.
It will be used in Section 2 to prove Theorem \ref{thm:main} on the
$\hbar$-adic valuation of universal $R$-matrices.

\medskip

Let $n,p$ be integers $\geq 1$.

\smallskip

{\bf - a - The Lie algebra $F_{(n,p)}$}

\smallskip

Let us denote by $\Fnp$ the free Lie algebra with generators $\widetilde{x}_{i,j}$,
where $(i,j)\in\{1,\dots,n\}\times\{1,\dots,p\}$.
Let us denote by $F_{(n,p)}$ the
quotient of
$\Fnp$ by the relations 
$$[\widetilde{x}_{i,j},\widetilde{x}_{i',j'}]=0\hbox{ when }i\not=i'\hbox{ and }
j\not=j'.$$
We denote by $x_{i,j}$ the image of $\widetilde{x}_{i,j}$ in $\Fnp$.
The Lie algebra $\Fnp$ is graded by $\Oplus\limits_{(i,j)
\in\{1,\dots,n\}\times\{1,\dots,p\}}\NM \epsilon_{(i,j)}$ where
$\deg(\widetilde{x}_{i,j})=\epsilon_{(i,j)}$.
We call this grading the {\it fine degree}.
$\Fnp$ has also a $\NM$-grading which we call the {\it total degree};
the total degree of $\widetilde{x}_{i,j}$ is $1$.
These gradings induce gradings on $F_{(n,p)}$.
\newline
For $A=\Fnp$ or $F_{(n,p)}$ and $\underline{k} \in
\oplus_{(i,j)} \NM \epsilon_{(i,j)}$ (resp., $k \in \NM$), we denote by $A_{\underline{k}}$
(resp., $A_k$) the fine degree $\underline{k}$ part (resp., total degree $k$ part)
of $A$.
We denote by $\widehat{A}$ the completion of $A$ with respect to the total degree.
So
$$\widehat{A}=\prod_{k\geq 0} A_k.$$

\smallskip

{\bf - b - The Campbell-Baker-Hausdorff series}

\stepcounter{subsection} \label{CBH}

\smallskip

If $N$ is an integer, let us denote by $\F_N$ the free Lie algebra with
generators $\widetilde{x}_1,\dots,\widetilde{x}_N$ and by $\widehat{\F}_N$ its
completion with respect to the total degree. There exists a series
$$\widetilde{x}_1\star \cdots \star \widetilde{x}_N \in \widehat{\F}_N,$$
such that the identity
$$\exp(\widetilde{x}_1\star \cdots \star \widetilde{x}_N )=\exp(\widetilde{x}_1)\cdots
\exp(\widetilde{x}_N)$$
holds in the completion $\widehat{\mathop{\hbox{\rm Freealg}}\nolimits}_N$ of
the free associative algebra 
$\mathop{\hbox{\rm Freealg}}\nolimits_N$ with generators
$\widetilde{x}_1,\dots,\widetilde{x}_N$, with respect to the total degree (see
\cite{Bo}).

\medskip

Let $\g$ be a $\NM$-graded Lie algebra, so $\g=\Oplus_{n \geq 0} \g_n$ and
let $\widehat{\g}$ be its completion; so $\widehat{\g}=\prod\limits_{n \geq 0} \g_n$.
If $x_1,\dots,x_n$ are elements of $\g$, with valuation $>0$, then there is a unique
continuous Lie
algebra
morphism $\pi$~:
$\widehat{\F}_N \to \widehat{\g}$, taking each $\widetilde{x}_i$ to $x_i$.
We then define
$$x_1\star \cdots \star x_N:=\pi(\widetilde{x}_1 \star \cdots \star \widetilde{x}_N).$$

\bigskip

{\bf - c - The elements $\delta_{(n,p)}$ of $\widehat{F}_{(n,p)}$}

\stepcounter{subsection} 

\bigskip

\begin{definition} We set
$$\delta_{(n,p)}=\sum_{k=0}^n \hskip0.5cm \sum_{l=0}^p 
\sum_{{\begin{array}{cl}
& {\scriptstyle{1\leq i_1<\cdots<i_k\leq n,}}\cr
& {\scriptstyle{1\leq j_1<\cdots<j_k\leq p}}
\end{array}}}{(-1)}^{n+p-k-l} \hskip4cm
$$
$$\hskip2cm (x_{i_1,j_1}\star
x_{i_1,j_2}\star \cdots \star x_{i_1,j_l})\star
(x_{i_2,j_l}\star\cdots \star x_{i_2,j_l})\star \cdots \star
(x_{i_k,j_1}\star\cdots \star x_{i_k,j_l}).$$
According to Section 1-b, the element $\delta_{(n,p)}$ belongs to
$\widehat{F}_{(n,p)}$.
\end{definition}

\medskip

\begin{theorem} \label{valuation}
The valuation of $\delta_{(n,p)}$ for the total degree of
$\widehat{F}_{(n,p)}$ is $\geq n+p-1$.
\end{theorem}

\smallskip

\dem
For $\underline{k} \in \Oplus_{(i,j)\in \{1,\dots,n\}\times\{1,\dots,p\}}
\NM \epsilon_{(i,j)} $,
we denote by $\delta_{(n,p),\underline{k}}$ the fine degree $\underline{k}$
component of $\delta_{(n,p)}$. Define the support of $\underline{k}$, 
$\supp(\underline{k})$ as the set of all pairs $(i,j)$ such that ${(i,j)}\not=0$.
If $S$ is a subset of $\{1,\dots,n\}\times\{1,\dots,p\}$, we call the
{\it $i$-th column of $S$} the intersection of $S$ with the $i$-th column
$\{i\}\times\{1,\dots,p\}$, and the 
{\it $j$-th line of $S$} its intersection with the $j$-th line
$\{1,\dots,n\}\times \{j\}.$
If $S_1$ is a subset of $\{1,\dots,n\}$, we denote by $\bar{S}_1$ its
complement $\{1,\dots,n\}\setminus S_1$;
if $S_2$ is a subset of $\{1,\dots,p\}$, we
denote by $\bar{S}_2$ its complement $\{1,\dots,p\}\setminus S_2$. We will
show:

\begin{proposition}
\label{1}
For each $\underline{k} \in \oplus_{(i,j)} ~\NM \epsilon_{(i,j)}$, we have
$\delta_{(n,p),\underline{k}}=0$ unless $S=\supp(\underline{k})$ satisfies the
following conditions:

\noindent \hbox{~}\hskip0.2cm(1)~each column of $S$ is nonempty;

\noindent \hbox{~}\hskip0.2cm(2)~each line of $S$ is nonempty;

\noindent \hbox{~}\hskip0.2cm(3)~if $S_1\subset\{1,\dots,n\}$ and 
$S_2\subset\{1,\dots,p\}$ are proper subsets
(i.e., they are non-empty
as well as their complements)
then $S\not\subset(S_1\times S_2)\cup
(\bar{S}_1 \times \bar{S}_2).$ 
\end{proposition}

\begin{proposition}
\label{combilemma}
If $S$ is a subset of $\{1,\dots,n\}\times\{1,\dots,p\}$ satisfying
conditions $(1)$, $(2)$ and $(3)$ of Proposition \ref{1}, then $\Card(S) \geq n+p-1$.
\end{proposition}

\noindent Theorem \ref{valuation} now follows from these propositions 
and the fact that the
total degree of an element of $\Fnp$ of degree $\underline{k}$ is
$\geq \Card(\supp(\underline{k}))$.
\findem

\vskip20pt

{\bf - d - {\sc Proof of Proposition \ref{1}.}}

\stepcounter{subsection}

\medskip 

For $i \in \{1,\dots,n\}$, 
there is a unique continuous Lie algebra morphism 
$$\lambda_i~:~\widehat{F}_{n,p} \to \widehat{F}_{n-1,p}$$
$$\hskip3.4cm x_{i',j'} \mapsto \left \{
\begin{array}{cl}
x_{i',j'}& \hbox{if }i'<i\cr
0& \hbox{if }i'=i\cr
x_{i'-1,j'}& \hbox{if }i'>i. \cr
\end{array}
\right.
$$
Similarly, for $j \in \{1,\dots,p\}$, 
there is a unique continuous Lie algebra morphism
$$\rho_j~:~ \widehat{F}_{n,p} \to \widehat{F}_{n,p-1}$$
$$\hskip3.4cm x_{i',j'} \mapsto \left \{
\begin{array}{cl}
x_{i',j'}& \hbox{if }j'<j\cr
0& \hbox{if }j'=j\cr
x_{i',j'-1}& \hbox{if }j'>j. \cr
\end{array}
\right.
$$
The proof will essentially consist in the following Lemmas \ref{lemma:1},
\ref{lemma:2}, and \ref{lemma:3}.
\begin{lemma} \label{lemma:1}
For each $i \in \{1,\dots,n\}$, we have $\lambda_i(\delta_{(n,p)})=0$.

\noindent
For each $j\in \{1,\dots,p\}$, we have $\rho_j(\delta_{(n,p)})=0$.
\end{lemma}

\smallskip

\begin{lemma} \label{lemma:2}
If $\alpha \in \widehat{F}_{(n,p)}$ is such that $\lambda_i(\alpha)=0$,
then the homogeneous component $\alpha_{\underline{k}}$ of $\alpha$
satisfies
$\alpha_{\underline{k}}=0$, unless
the $i$-th line of $\supp(\underline{k})$ is nonempty.

\noindent
In the same way, if $\rho_j(\alpha)=0$, then
$\alpha_{\underline{k}}=0$ unless
the $j$-th line of $\supp(\underline{k})$ is nonempty.
\end{lemma}

\smallskip

\noindent After Lemmas \ref{lemma:1} and \ref{lemma:2} are proved, we know that
$\delta_{(n,p),\underline{k}}=0$ unless each line
and each column of $S=\supp(\underline{k})$ is
nonempty. This proves that $S$ should
satisfy conditions $(1)$ and $(2)$. The fact that $S$ satisfies condition $(3)$
will follow from:

\begin{lemma} \label{lemma:3}
Let $\alpha \in \widehat{F}_{(n,p)}$ be homogeneous of
degree $\underline{k}$.
Assume that each line and each column of $\supp(\underline{k})$
is nonempty.
Then $\alpha=0$ unless $\underline{k}$ satisfies conditions
$(3)$ of Proposition \ref{1}.
\end{lemma}

\medskip

\noindent So Lemmas \ref{lemma:1},
\ref{lemma:2}, and \ref{lemma:3} imply Proposition \ref{1}. 
\findem

\medskip

We now prove Lemmas \ref{lemma:1},
\ref{lemma:2}, and \ref{lemma:3}:

\medskip

\noindent{\sc Proof of Lemma \ref{lemma:1}.}
If $I \subset \nn$ and
$J\subset \pp$,
let us set $k=\Card(I)$, $l=\Card(J)$, and
$I=\{i_1,\dots,i_k\}$, $J=\{j_1,\dots,j_l\}$ where
$1\leq i_1 < \cdots < i_k \leq n$ and
$1\leq j_1 < \cdots < j_l \leq p$. Let us then set
$$\delta_{(n,p)}^{I,J}={(-1)}^{n+p-k-l}
x_{i_1,j_1} \star x_{i_1,j_2} \star \cdots \star x_{i_k,j_l}.$$
Then we have:
$$\delta_{(n,p)}=\sum_{\begin{array}{rl} &{\scriptstyle I \subset \nn ,}\cr &
{\scriptstyle J \subset \pp}
\end{array}} \delta_{(n,p)}^{I,J}.$$
Let $i \in \nn$ and let us prove that $\lambda_i(\delta_{(n,p)})=0$:
we write
$$\delta_{(n,p)}=\sum_{\begin{array}{rl} &
{\scriptstyle I'\subset \nn\setminus\{i\},}\cr
&{\scriptstyle J \subset \pp} \end{array}} \delta_{(n,p)}^{I',J}+
\delta_{(n,p)}^{I'\cup\{i\},J},$$
and we now show that for any $(I',J)$,
$$\lambda_i\left(\delta_{(n,p)}^{I',J}+
\delta_{(n,p)}^{I'\cup\{i\},J}\right)=0.$$
This follows from the following facts:
\begin{itemize}
\item if $\varphi$~: $\widehat{\g} \to \widehat{\g}'$ is a morphism of complete graded
Lie algebras,
and $x_1,\dots,x_N$ are elements of $\widehat{\g}'$ of positive valuation, then
$$\varphi(x_1)\star\cdots\star\varphi(x_N)=\varphi(x_1\star\cdots\star x_n);$$
\item we have the identity $x_1 \star \cdots \star x_k \star 0 \star x_{k+1}
\star \cdots \star x_N= x_1 \star \cdots \star x_N$.
\end{itemize}
The proof of $\rho_j(\delta_{(n,p)})=0$ is similar: we write
$$\delta_{(n,p)}=\sum_{\begin{array}{rl} &
{\scriptstyle I\subset \nn,}\cr
&{\scriptstyle J' \subset \pp \setminus \{j\}}\end{array}} \delta_{(n,p)}^{I,J'\cup\{j\}}
+ \delta_{(n,p)}^{I,J'},$$
and then show that 
$$\rho_j\left(\delta_{(n,p)}^{I,J'\cup\{j\}}
+ \delta_{(n,p)}^{I,J'}\right)=0$$
for any pair $(I,J')$ using the same arguments.
\findem

\bigskip

\noindent{\sc Proof of Lemma \ref{lemma:2}.}
It will be enough to prove Lemma \ref{lemma:2} when $i=n$.
Let $\alpha$ be an element of $F_{n,p}$.
Set
$$\alpha=\sum_{\uk \in \Oplus_{(i,j)}~\NM \epsilon_{(i,j)}}
\alpha_{\uk},$$
where $\alpha_{\uk}$ has degree $\uk$. Assume that $\lambda_n(\alpha)=0$,
we want to show that $\alpha_{\uk}=0$ unless the $n$-th column of
$\supp(\uk)$ is nonempty. We will now write:
$$\alpha=\sum_{
\uk~|~n{\scriptstyle\rm{\hbox{-}th~ column(supp}}(\uk))\not=\emptyset}\alpha_{\uk}+
\sum_{\uk~|~n{\scriptstyle\rm{\hbox{-}th~ column(supp}}(\uk))=\emptyset}\alpha_{\uk}.$$
Now $\lambda_n(\alpha_{\uk})=0$ as soon as the $n$-th column of $\supp(\uk)$
is nonempty, so we get
\begin{equation}
\sum_{\uk~|~n{\scriptstyle\rm{\hbox{-}th~column(supp}}(\uk))=\emptyset}\lambda_n(\alpha_{\uk})=0
\hskip0.5cm (\hbox{identity in }F_{(n-1,p)}).
\label{I}
\end{equation}
There is a unique Lie algebra morphism
$\iota$~: $F_{n-1,p} \to F_{n,p}$, taking each $x_{i,j}$ to $x_{i,j}$.
We have 
$$\lambda_n \circ \iota = \Id_{F_{n-1,p}},$$ 
therefore $\iota$ is
injective. Moreover, let $\FC$ be the Lie subalgebra of $F_{n,p}$
generated by the $x_{i,j}$, $i\in \nn$, $j \in \pp$.
$\FC$ is the image of the Lie subalgebra $\widetilde{\FC} \subset \Fnp$
generated by the $\widetilde{x}_{i,j}$, $i\in \{1,\dots,n-1\}$,
$j\in \pp$ under the canonical projection $\Fnp \to F_{(n,p)}$.
Denote by $\widetilde{\iota}$~: 
$\mathop{\hbox{\rm Free}}\nolimits_{(n-1,p)}\to \Fnp$ the morphism taking each
$\widetilde{x}_{i,j}$ to $\widetilde{x}_{i,j}$, then the diagram
$$\begin{matrix}
{\mathop{\hbox{\rm Free}}\nolimits_{(n-1,p)}}
&\mathop{\longrightarrow}\limits^{\widetilde{\iota}}
& {\widetilde{\FC}}
& \hookrightarrow
&{\Fnp}\cr
\downarrow&&\downarrow&&\downarrow\cr
{F_{(n-1,p)}}
&\mathop{\longrightarrow}\limits^\iota  
&{\FC}
&\hookrightarrow 
& { F_{(n,p)}}
\end{matrix}$$
commutes. Since the vertical arrows are onto and
$\widetilde{\iota}$~: $\mathop{\hbox{\rm Free}}\nolimits_{(n-1,p)}\to
\widetilde{\FC}$ is onto, $\iota$~: $\mathop{\hbox{\rm Free}}\nolimits_{(n-1,p)}
\to \FC$ is onto. So $\iota$~: $\mathop{\hbox{\rm Free}}\nolimits_{(n-1,p)}
\to \FC$ is an isomorphism. Now if the $n$-th column of
$\supp(\uk)$ is empty, we have $\alpha_{\uk} \in \FC$.
Apply $\iota$ to identity (\ref{I}), we get
$$\sum_{\uk~|~n{\scriptstyle \rm{\hbox{-}th~ column(supp}}(\uk))=\emptyset}\alpha_{\uk}=0
\hskip0.5cm (\hbox{identity in }F_{(n-1,p)}).$$
Separating homogeneous components, we get $\alpha_{\uk}=0$ for each $\uk$, such
that
$n$-th column$(\supp(\uk))=\emptyset$.
\findem

\bigskip

\noindent{\sc Proof of Lemma \ref{lemma:3}.}
Let $\alpha \in F_{(n,p)}$ be of fine degree $\uk$, and assume that $\uk$
satisfies $(1)$, $(2)$ but not $(3)$. So we have
proper subsets $S_1 \subset \nn$ and
$S_2 \subset \pp$ such that
$$\supp(\uk) \subset (S_1\times S_2) \cup (\bar{S}_1\times \bar{S}_2).$$
The rectangles $R=S_1\times S_2$ and
$R'=\bar{S}_1\times \bar{S}_2$ are disjoint.
Moreover, condition $(1)$ (or condition $(2)$) implies
that $\supp(\uk)\cap R$ and
$\supp(\uk)\cap R'$ are both nonempty.
Let $(k_{(i,j)})_{(i,j)}$ be the component of $\uk$ in the
basis 
$(\epsilon_{(i,j)})_{(i,j)\in \nn \times \pp}$. Then
$k_{(i,j)}\not=0$ if and only if $(i,j)\in \supp(\uk)$, so
$$\uk = \sum_{(i,j)\in {\scriptstyle \rm{supp}}(\uk)}k_{(i,j)}\epsilon_{(i,j)}.$$
Now there exists $\beta \in \Fnp$ of fine degree $\uk$, whose image under
$\Fnp \to F_{(n,p)}$ is $\alpha$. Moreover, standard results on free Lie algebra
(see \cite{Bo}) imply that $\beta$ has the following form:
set $|\uk|=\sum_{(i,j)\in \nn \times \pp}k_{(i,j)}$, and say that a map
$\mu$~: $\{1,\dots,|\uk|\}\to \nn \times \pp$ is a $\uk$-map if for
any $(i,j)\in \nn \times \pp$, $\Card(\mu^{-1}(i,j))=k_{(i,j)}$.
If  $\mu$ is a $\uk$-map, let us define $\widetilde{x}_\mu$ as the iterated Lie
bracket
$$\widetilde{x}_\mu=\left[\dots
\left[\left[\widetilde{x}_{\mu(1)},\widetilde{x}_{\mu(2)}\right],
\widetilde{x}_{\mu(3)}\right],
\dots,\widetilde{x}_{\mu(|\uk|)}\right].$$
Then $\beta$ is a linear combination
$$\beta=\sum_{\mu \in \{\uk{\scriptstyle \rm{\hbox{-}maps}}\}}\beta_\mu\widetilde{x}_\mu,$$
where the $\beta_\mu$ are scalars. Define $x_\mu$ as the image of
$\widetilde{x}_\mu$ under the map $\Fnp \to F_{(n,p)}$.
Then $\alpha=\sum_{\mu \in \{\uk{\scriptstyle \rm{\hbox{-}maps}}\}}\beta_\mu x_\mu$. On the
other hand, $x_\mu=0$ if $\mu$ is any $\uk$-map. Indeed
$$x_\mu=\left[\dots
\left[\left[{x}_{\mu(1)},{x}_{\mu(2)}\right],{x}_{\mu(3)}\right],
\dots,{x}_{\mu(|\uk|)}\right].$$
Then, since $R\cap R' = \emptyset$, $\supp(\uk)\cap R\not=\emptyset$
and $\supp(\uk)\cap R'\not=\emptyset$, there
exists an integer $\nu \in \{1,\dots,|\uk|-1\}$, such that
one of the following possibilities occurs:
\begin{itemize}
\item either $\mu(1),\dots,\mu(\nu) \in R$ and $\mu(\nu+1) \in R'$
\item or $\mu(1),\dots,\mu(\nu) \in R'$ and $\mu(\nu+1) \in R$.
\end{itemize}
Now the bracket $[x_\alpha,x_\beta]$ vanishes when $\alpha \in R$ and
$\beta \in R'$. It follows that in both cases, the bracket
$$\left[\dots
\left[\left[{x}_{\mu(1)},{x}_{\mu(2)}\right],{x}_{\mu(3)}\right],
\dots,{x}_{\mu(\nu+1)}\right]$$
vanishes. Therefore $x_\mu=0$, which implies $\alpha=0$.
\findem

\bigskip

{\bf - e - {\sc Proof of Proposition \ref{combilemma}.}}

\stepcounter{subsection}

\medskip

The argument will be an induction over $n+p$. Assume that the proposition is
proved when 
$n+p \leq N$ and let us prove it when $n+p=N+1$.
Let $(n,p)$ be such that $n+p=N+1$ and let $S$ be as in the proposition.
We may assume $n \geq p$. If each column of $S$ has $\geq 2$ elements, then
$$\Card(S) \geq n.2 \geq n+p-1.$$
Each column of $S$ has $\geq 1$ element by assumption $(1)$, so we can
assume that
$S$ has a column with exactly $1$ element.
We may assume that this element is $(n,p)$.
Let us now set $S'=S\cap(\nn \times \pp)$ and prove that $S'$ satisfies the
analogues (1'), (2') and (3') of the assumptions of the proposition, 
where $(n,p)$ is replaced by $(n-1,p)$:

\noindent \hbox{~}\hskip0.2cm(1')~ for $i=1,\dots,n-1$, the $i$-th lines of $S$ and of 
$S'$ coincide so the $i$-th
line of $S'$ is nonempty;

\noindent \hbox{~}\hskip0.2cm(2') for $j=1,\dots,p-1$, the $j$-th lines of $S$
and of $S'$ coincide, so the 
$j$-th line of $S'$ is nonempty. Let us prove that the $p$-th line of $S'$ is
also nonempty. If not, we would have 
$(n$-th column of $S)\cup(p$-th line of $S)=\{(n,p)\}$,
contradicting assumption 
$(3)$ on $S$ with $S_1=\{1,\dots,n-1\}$ and 
$S_2=\{1,\dots,p-1\}$. Therefore all the
lines
of $S'$ are nonempty;

\noindent \hbox{~}\hskip0.2cm(3')~ assume that there exist proper subsets 
$S_1'\subset
\{1,\dots,n-1\},~S_2\subset \pp$, (i.e., nonempty and with nonempty
complement), such that 
$$S'\subset(S_1'\times S_2)\cup\left((\{1,\dots,n-1\}\setminus S_1')\times
(\pp\setminus S_2)\right).$$ 
Then if 
$p\in S_2$, let us set 
$S_1=S_1  \cup\{n\}$, and if 
$p\in \pp\setminus S_2$, let us set 
$S_1=S_1'$.
Then in each case, since 
$S=S'\cup\{(n,p)\}$,
we get
$$S \subset (S_1\times S_2)\cup
\left((\nn \setminus S_1)\times (\pp\setminus S_2)\right),$$
thus contradiction the fact that
$S$ satisfies assumption $(3)$.
Therefore $S'$ satisfies assumption (3').

\medskip

Now $S'$ satisfies the assumptions of the proposition, with $(n,p)$ replaced
with
$(n-1,p)$.
Since $(n-1)+p=N$, we may apply the induction hypothesis. So
$\Card(S')\geq (n-1)+p-1$.
Now $\Card(S)=\Card(S')+1$, so $\Card(S)\geq n+p-1$. This proves the induction
step.
\hbox{~}~\findem

\vskip20pt

\centerline {\bf \S\; 2 \ $\hbar$-adic valuation properties of universal
$R$-matrices }

\stepcounter{section}\label{part2}

\vskip20pt

In this section, we first identify a tensor product of QFSH 
(quantized formal series Hopf) algebras as a
subspace of the tensor product of the
corresponding QUE algebras (Proposition \ref{pdt:QFSH}).

\medskip

{\bf - a - {Tensor product of QFSH algebras}}

\stepcounter{subsection}

\medskip

Let $H$ be a QUE algebra.
The corresponding QFSH algebra is:
$$H'=\{x\in H|~\forall n \in \NM,~ \delta_H^{(n)}(x) \in \hbar^n 
H^{\widehat{\otimes} n}\}$$ 
(where $\delta_H^{(n)}=\left(\Id - \eta \circ \varepsilon\right)^{\otimes
n}\circ \Delta_H^{(n)}$).
We express this condition as follows: set
$$\delta_H^{(\tot)}=\prod_{n \geq 0}\delta_H^{(n)}~:~H\to \prod_{n\geq
0}H^{\widehat{\otimes} n}.$$
Then 
$$H'=\left(\delta_H^{(\tot)}\right)^{-1}\left(\prod_{n\geq
0}\hbar^n H^{\widehat{\otimes} n}\right).$$
In the same way, if $K$ is another QUE algebra, we define
$$K'=\left(\delta_K^{(\tot)}\right)^{-1}\left(\prod_{n\geq
0}\hbar^n K^{\widehat{\otimes} n}\right).$$
Then $\delta_H^{(\tot)}\otimes
\delta_K^{(\tot)}$ is a linear map $H\widehat{\otimes}K \to
\prod\limits_{n,p\geq0}
H^{\widehat{\otimes} n}\widehat{\otimes}K^{\widehat{\otimes} p}.$

\begin{proposition}\label{pdt:QFSH}
The tensor product of QFSH algebras $H'$ and $K'$ is
$$H' \bar{\otimes}\ K'=\left(\delta_H^{(\tot)}\otimes
\delta_K^{(\tot)}\right)^{-1}\left(\prod_{n,p\geq
0}\hbar^{n+p} H^{\widehat{\otimes} n}\widehat{\otimes}K^{\widehat{\otimes} p}\right).$$
\end{proposition}

\dem
We will need the following lemma:
\begin{lemma}
Let $V_1,~V_2,~T_1,~T_2$ be vector spaces and let 
$\delta_1$~:
$V_1 \hookrightarrow T_1$,
$\delta_2$~: $V_2\hookrightarrow T_2$ be injections. Let $S_1,~S_2$ be vector
subspaces of $T_1,~T_2$, and set $U_i=\delta_i^{-1}(S_i)$.
Then
$$U_1 \otimes U_2=(\delta_1 \otimes \delta_2)^{-1}(S_1 \otimes S_2).$$
\label{lemmm}
\end{lemma}

The proposition then follows because the maps $\delta_K^{(\tot)}$
and $\delta_H^{(\tot)}$ are injective.
Indeed, $x$ can be obtained from $\delta_H^{(\tot)}(x)$ by the formula:
$$x=\left(\delta_H^{(\tot)}(x)\right)_1+ \eta \left(
\left(\delta_H^{(\tot)}(x)\right)_0\right).$$
\hbox{~}\findem

\medskip

\noindent{\sc Proof of Lemma \ref{lemmm}.} Let $W_1=\Im(\delta_1)$,
$W_2=\Im(\delta_2)$. Then $W_1 \otimes W_2=\Im(\delta_1 \otimes \delta_2)$. So
$$W_1 \subset T_1,\hskip0.5cm
W_2 \subset T_2,\hskip0.5cm
W_1 \otimes W_2 \subset T_1 \otimes T_2.$$ 
We should prove that
if $W_1 \subset T_1$, $W_2 \subset T_2$,
then the subspace $\left(W_1 \cap S_1\right) \otimes 
\left(W_2 \cap S_2 \right)$
$\subset W_1 \otimes W_2$ is equal to 
$(W_1 \otimes W_2) \cap (S_1 \otimes S_2)$.
Identifying the spaces $V_i$ with their images $W_i=\Im(\delta_i)$ through the
maps $\delta_i$, we should prove the following statement:

\begin{lemma}
Let $T_1,~T_2$ be vector spaces, and let for each $i=1,~2$, $W_i$ and $S_i$ be
vector subspaces of
$T_i$. Then
$$(W_1 \cap S_1)\otimes (W_2 \cap S_2)=(W_1 \otimes W_2) \cap (S_1 \otimes
S_2).$$
\label{lembis}
\end{lemma}

Indeed, $U_1 \otimes U_2$ can be identified with $(W_1 \cap S_1)\otimes (W_2
\cap S_2)$ through
$\delta_1 \otimes \delta_2$, whereas $(\delta_1 \otimes \delta_2)^{-1}(S_1
\otimes S_2)$ can be identified with $(W_1 \otimes W_2)\cap (S_1 \otimes S_2)$
through the same map.
\findem

\medskip

\noindent{\sc Proof of Lemma \ref{lembis}.} Introducing
supplementary vector spaces $\widetilde{W}_i$ of $W_i \cap S_i$ in $W_i$ and
$\widetilde{S}_i$ of $W_i \cap S_i$ in $S_i$, we get
$$\left \{
\begin{array}{cl}
W_i=& (W_i\cap S_i)\oplus \widetilde{W}_i\cr
S_i=& (W_1\cap S_i)\oplus \widetilde{S}_i
\end{array}
\right.
$$
$W_i+S_i=(W_i\cap S_i)\oplus \widetilde{W}_i \oplus \widetilde{S}_1$ so
$(W_1+S_1)\otimes(W_2+S_2)$ is the direct sum of nine summands.
The subspaces $W_1 \otimes W_2$ and $S_1 \otimes S_2$ are the sums
$$\left \{
\begin{array}{rl}
W_1\otimes W_2=& \big((W_1\cap S_1)\otimes (W_2\cap S_2)\big)\oplus 
\big((W_1\cap S_1)\otimes \widetilde{W}_2 \big)\cr
&\oplus 
\big(\widetilde{W}_1 \otimes (W_2\cap S_2)\big) \oplus
\big(\widetilde{W}_1 \otimes \widetilde{W}_2\big)\cr
S_1\otimes S_2=& \big((W_1\cap S_1)\otimes (W_2\cap S_2)\big)\oplus 
\big((W_1\cap S_1)\otimes \widetilde{S}_2 \big)\cr
&\oplus ~
\big(\widetilde{S}_1 \otimes (W_2\cap S_2)\big) \oplus
\big(\widetilde{S}_1 \otimes \widetilde{S}_2 \big),
\end{array}
\right.
$$
the intersection of which is $(W_1\cap S_1)\otimes (W_2\cap S_2)$.
\findem

\medskip

{\bf - b - {\sc Proof of Theorem \ref{thm:main}.}}

\stepcounter{subsection}
\label{section:proof:main}

\medskip

Let $(U_\hbar(\g),R)$ be a quasitriangular QUE algebra.
Set
$\rho=\hbar\log(R)$. Then $\rho \in \hbar^2 U_\hbar(\g)^{\widehat{\otimes}2}$.
For $x,~y \in U_\hbar(\g)^{\widehat{\otimes}n}[\hbar^{-1}]$, let us set
$$\{x,y\}_\hbar={{1}\over{\hbar}}[x,y].$$
Then $\{-,-\}$ restricts to a Lie bracket
$$\{-,-\}~:~\ve^2\left(\hbar^2U_\hbar(\g)^{\widehat{\otimes}k}\right)
\to \hbar^2U_\hbar(\g)^{\widehat{\otimes}k}.$$
Its image is actually contained in $\hbar^3U_\hbar(\g)^{\widehat{\otimes}k}$.
According to Proposition \ref{pdt:QFSH}, the statement
$\rho \in \OC_\hbar(G^*)^{\bar{\otimes}2}$ is equivalent to 
\begin{equation}\left(\delta^{(n)}_{U_\hbar(\g)} \otimes \delta^{(p)}_{U_\hbar(\g)} \right)
(\rho) \in \hbar^{n+p} U_\hbar(\g)^{\widehat{\otimes}(n+p)} 
\label{S}
\end{equation}
for any $n,~p \geq 0$. Let us therefore prove this statement:
for $a,b \in \hbar^2 U_\hbar(\g)^{\widehat{\otimes}k}$, let us set
\begin{equation}
a \star_kb=a+b+{{1}\over{2}}\{a,b\}_\hbar +\cdots~;
\label{CBH1}
\end{equation}
the series (\ref{CBH1}) makes sense because
of the following fact: if $m \in \mathop{\hbox{\rm Free}}\nolimits_{2}$ is a Lie
monomial of degree $p$ in two free variables $A$, $B$, and
$m\left(\{-,-\}_\hbar,a,b\right)$ is the image of $m$ by the Lie algebra
morphism $\mathop{\hbox{\rm Free}}\nolimits_{2}\to
\hbar^2U_\hbar(\g)^{\widehat{\otimes}k}$ defined by $A \mapsto a$, $B\mapsto b$,
then
$$m\left(\{-,-\}_\hbar,a,b\right)\in \hbar^{p+1}U_\hbar(\g)^{\widehat{\otimes}k}.$$
Recall that $\rho \in \hbar^2U_\hbar(\g)^{\widehat{\otimes}3}$. The quasitriangular
identities then imply that $\rho$ satisfies
\begin{equation}
\label{ID1}
\left(\Delta_{U_\hbar(\g)}\otimes
\Id\right)(\rho)=\rho^{1,3}\star_3\rho^{2,3},~
\left(\Id\otimes\Delta_{U_\hbar(\g)}\right)(\rho)=\rho^{1,3}\star_3\rho^{1,2}.
\end{equation}
There is a unique Lie algebra morphism 
$$\phi_{(n,p)}~:~
\widehat{F}_{(n,p)}\to
\left(\hbar^2 U_\hbar(\g)^{\widehat{\otimes}(n+p)},\{-,-\}_\hbar\right),$$
$$\hskip-1.5cm x_{i,j} \mapsto \rho^{i,n+j}$$ 
for each 
$(i,j)\in \nn\times \pp$
(here $\rho^{i,n+j}$ is the image of $\rho$ by the map
$$U_\hbar(\g)^{\widehat{\otimes}2}\to U_\hbar(\g)^{\widehat{\otimes}(n+p)},$$
taking the first component of $U_\hbar(\g)^{\widehat{\otimes}2}$ to the $i$-th
component, and the second component to the $(n+j)$-th component).
Then the identities $(\ref{ID1})$ imply
\begin{equation}
\label{ID2}
\left(\delta_{U_\hbar(\g)}^{(n)}\otimes \delta_{U_\hbar(\g)}^{(p)}\right)
(\rho)=\phi_{(n,p)}\left(\delta_{(n,p)}\right).
\end{equation}
Let $\alpha \in \F_{(n,p)}$ be an element of total degree $d$.
Then we have seen that
$\phi_{(n,p)}(\alpha)\in\hbar^{d+1}U_\hbar(\g)^{\widehat{\otimes}(n+p)}$.
Since $\delta_{(n,p)}$ has valuation $\geq n+p-1$, we get 
$\phi_{(n,p)}(\delta_{(n,p)})\in \hbar^{n+p}
U_\hbar(\g)^{\widehat{\otimes}(n+p)}$\!.
According to $(\ref{ID2})$, this implies identity $(\ref{S})$.
So we get $\rho \in \OC_\hbar(G^*)^{\bar{\otimes}2}$.
Since we have $(\varepsilon\otimes \Id)(R)=(\Id \otimes \varepsilon)(R)=1$,
we get $(\varepsilon \otimes \Id)(\rho)=(\Id \otimes \varepsilon)(\rho)=0$, therefore 
$\rho\in (\m_\hbar)^{\bar{\otimes}2}$.
\findem

\vskip20pt

\centerline {\bf \S\; 3 \ Proof of Corollary \ref{cor:GH}.}

\stepcounter{section}

\vskip20pt

The space $U_\hbar(\g)^{\widehat{\otimes}2}$ is a quantization
of $\g \times \g$ and
$\OC_\hbar(G^*)^{\bar{\otimes}2}$ is the corresponding QFSH algebra.
Set $\d=\g \times \g$, let $D^*$ be the corresponding formal group, and
set $U_\hbar(\d)=U_\hbar(\g)^{\widehat{\otimes}2}$,
$\OC_\hbar(D^*)=\OC_\hbar(G^*)^{\bar{\otimes}2}$. Denote by $\m_\hbar(D^*)$
the augmentation ideal of
$\OC_\hbar(D^*)$. Then $\m_\hbar \bar{\otimes}1$ and
$1 \bar{\otimes} \m_\hbar$ are subspaces of $\m_\hbar(D^*)$,
so $\rho \in \m_\hbar(D^*)^2$. Corollary \ref{cor:GH} will follow from the
following proposition:

\begin{proposition}
Let $\d$ be an arbitrary finite-dimensional Lie bialgebra, let 
$U_\hbar(\d)$ be a quantization of $\d$, let
$\OC_\hbar(D^*)$ be the QFSH subalgebra of $U_\hbar(\d)$, and
let
$\m_\hbar(D^*)$ be
the augmentation ideal of
$\OC_\hbar(D^*)$. Let $\rho$ be an arbitrary element of $\m_\hbar(D^*)^2$.
 
\noindent \hbox{~}\hskip0.2cm(1) the operation $\{-,-\}_\hbar$~: $\ve^2\left(
U_\hbar(\d)[\hbar^{-1}]\right) \to U_\hbar(\d)[\hbar^{-1}]$ restricts to
$$\{-,-\}_\hbar~:~ \ve^2\left(\OC_\hbar(D^*)\right) \to \OC_\hbar(D^*);$$

\noindent \hbox{~}\hskip0.2cm(2) for any $k,~l \geq 0$, we have
$$\left\{\m_\hbar(D^*)^k,\m_\hbar(D^*)^l\right\}_\hbar \subset
\m_\hbar(D^*)^{k+l-1};$$

\noindent \hbox{~}\hskip0.2cm(3) 
$\OC_\hbar(D^*)=\lim\limits_{\lar_k}\left(\OC_\hbar(D^*)/\big(
\m_\hbar(D^*)+\hbar \OC_\hbar(D^*)\big)^k\right)$~; 
let us set 
$$\ad_\hbar(\rho)(x):=\{\rho,x\}_\hbar,$$ 
then the exponential
$\exp\left(\ad_\hbar(\rho)\right)$ is a well-defined continuous algebra
automorphism of $\OC_\hbar(D^*)$;
 
\noindent \hbox{~}\hskip0.2cm(4) We have $\rho \in \hbar^2 U_\hbar(\d)$, so
$\exp\left({{\rho}\over{\hbar}}\right)$ is a well-defined element of
$1+\hbar U_\hbar(\d)$. Let 
$\Ad\left(\exp\left({{\rho}\over{\hbar}}\right)\right)$ be the inner automorphism 
$x \mapsto \exp\left({{\rho}\over{\hbar}}\right)x
\exp\left({{\rho}\over{\hbar}}\right)^{-1}$ of $U_\hbar(\d)$.
Then $\Ad\left(\exp\left({{\rho}\over{\hbar}}\right)\right)$ restricts to an
automorphism of $\OC_\hbar(D^*)$, which coincides with
$\exp(\ad_\hbar(\rho))$.
\end{proposition}

\medskip

\dem

\noindent \hbox{~}\hskip0.2cm(1) $\OC_\hbar(D^*)$ is a subalgebra of 
$U_\hbar(\d)$, and its reduction modulo
$\hbar$ is commutative, so for $x,~y \in \OC_\hbar(D^*)$,
$[x,y] \in \hbar \OC_\hbar(D^*)$, i.e., 
$\{x,y\}_\hbar\in\OC_\hbar(D^*)$.

\noindent \hbox{~}\hskip0.2cm(2) Let $\varepsilon$~: 
$\OC_\hbar(D^*)\to k[[\hbar]]$ be the
augmentation
map.
If $x,~y \in \OC_\hbar(D^*)$, then $\varepsilon([x,y])=0$. So
$\hbar\varepsilon(\{x,y\}_\hbar)=0$,
i.e., $\varepsilon(\{x,y\}_\hbar)=0$, so $\{x,y\}_\hbar\in \m_\hbar(D^*)$.
Therefore
$\{\m_\hbar(D^*),\m_\hbar(D^*)\} \subset \m_\hbar(D^*)$.
Applying the Leibniz rule, we get (2).

\noindent \hbox{~}\hskip0.2cm(3) The first part is contained in \cite{Ga}. 
The second part follows from
(2): $\rho \in \m_\hbar(D^*)^2$, therefore
$\ad_\hbar(\rho)\left((\m_\hbar(D^*)^k\right)\subset \m_\hbar(D^*)^{k+1}$ for
any $k$.

\noindent \hbox{~}\hskip0.2cm(4) If $x \in \OC_\hbar(D^*)$, 
then $x-\varepsilon(x) \in
\hbar U_\hbar(\d)$. If, in addition, $x\in \m_\hbar(D^*)$, then 
$\varepsilon(x)=0$ so $x \in \hbar U_\hbar(\d)$.
So $\m_\hbar(D^*) \subset \hbar U_\hbar(\d)$ and
$\m_\hbar(D^*)^2 \subset \hbar^2 U_\hbar(\d)$.
Let us now show that the diagram
\begin{equation}
\label{comm:diag}
\begin{matrix}
U_\hbar(\d)
& \To^{{\scriptstyle\rm{Ad}}\left(\exp\left({{\rho}\over{\hbar}}\right)\right)}
& U_\hbar(\d)\cr
\uparrow&&\uparrow\cr
\OC_\hbar(D^*)
&\mathop{\longrightarrow}\limits^{\exp\left({\scriptstyle\rm{ad}}_\hbar(\rho)\right)}
&\OC_\hbar(D^*)
\end{matrix}
\end{equation}
commutes. This means that we have the identity (in $U_\hbar(\d)$):
$$\exp\left(\ad_\hbar(\rho)\right)(x)=\exp\left(
{{\rho}\over{\hbar}}\right)x
\exp\left(
{{\rho}\over{\hbar}}\right)^{-1},$$
for any $x \in \OC_\hbar(D^*)$. To show this identity, let us introduce
a parameter $t$ and show the same identity, where $\rho$ is replaced by $t\rho$.
This means that we have the identity
$$\exp\left(\ad_\hbar(t\rho)\right)(x)=
\exp\left({{t\rho}\over{\hbar}}\right)x
\exp\left(
{{t\rho}\over{\hbar}}\right)^{-1}.$$
This last identity follows from the fact that its two sides
coincide when $t=0$, and that they both satisfy the differential equation
$${{{\rm{d}}}\over{{\rm{d}}\ t}}x(t)=\ad_\hbar(\rho)(x(t)).$$
It follows that the diagram (\ref{comm:diag}) commutes. Therefore the map
$\Ad\left(\exp\left({{\rho}\over{\hbar}}\right)\right)$ restricts to an
automorphism of $\OC_\hbar(D^*)$, which coincides with
$\exp\left(\ad_\hbar(\rho)\right)$.
\findem

\vskip1.1truecm

\centerline { ACKNOWLEDGEMENTS }

\vskip18pt

We would like to thank F.\ Gavarini for discussions on the subject of this
work.
We also would like to thank P.\ Cartier for suggesting that Theorem \ref{thm:main}
is formally analogous to the oscillatory integrals formalism.

\end{document}